\documentclass{conm-p-l}

\usepackage{amssymb}
\usepackage{enumerate}
\usepackage{bbm}

\newcommand{\excise}[1]{}
\newcommand{\comment}[1]{{$\star$\sf\textbf{#1}$\star$}}

\newtheorem{thm}{Theorem}[section]
\newtheorem{prop}[thm]{Proposition}
\newtheorem{lemma}[thm]{Lemma}
\newtheorem{cor}[thm]{Corollary}

\theoremstyle{definition}
\newtheorem{example}[thm]{Example}
\newtheorem{defn}[thm]{Definition}

\newtheorem{remark}[thm]{Remark}

\newcommand\0{\mathbf{0}}
\newcommand\1{\mathbbm{1}}
\newcommand\<{\langle}
\newcommand\CC{\mathbb C}
\newcommand\HH{{\widetilde H}{}{}}
\newcommand\NN{\mathbb N}
\newcommand\RR{\mathbb R}
\newcommand\ZZ{\mathbb Z}
\newcommand\bb{\mathbf b}
\newcommand\kk{\Bbbk}
\newcommand\mm{\mathfrak m}
\newcommand\uu{\mathbf u}
\newcommand\xx{\mathbf x}

\newcommand\CD{\text{\v C}\Delta}
\newcommand\ED{\mathrm{E}\Delta} 
\newcommand\KD{\mathrm{K}\Delta} 
 
\newcommand\oF{{\hspace{.3ex}\ol{\hspace{-.3ex}F\hspace{-.05ex}}\hspace{.05ex}}}
\newcommand\vC{\check{\mathcal{C}}}
\newcommand\from{\leftarrow}
\newcommand\spot{{\raisebox{.25ex}{\tiny$\scriptscriptstyle\bullet$}}}
\newcommand\minus{\smallsetminus}
\newcommand\goesto{\rightsquigarrow}
\newcommand\nothing{\varnothing}
\newcommand\Mustata{Musta\c t\v a}

\renewcommand\>{\rangle}
\renewcommand\aa{\mathbf a}
\renewcommand\iff{\Leftrightarrow}
\renewcommand\implies{\Rightarrow}

\DeclareMathOperator\Ext{Ext}
\DeclareMathOperator\Tor{Tor}

\DeclareMathOperator\qdeg{qdeg}
\DeclareMathOperator\supp{supp}
\DeclareMathOperator\tdeg{tdeg}

\newcommand\wt[1]{\widetilde{#1}{}}
\newcommand\ol[1]{{\overline{#1}}}
\newcommand\link[2]{\mathrm{link}\hspace{.2ex}(#1,#2)}
\newcommand\Star[2]{\mathrm{star}\hspace{.2ex}(#1,#2)}

\begin{document}

\mbox{}\vspace{15.73ex}
\title{Topological Cohen--Macaulay criteria for monomial ideals}
\author{Ezra Miller}
\address{School of Mathematics\\University of Minnesota\\Minneapolis, MN 55455}
\curraddr{Department of Mathematics\\Duke University\\Durham, NC 27708}
\email{ezra@math.duke.edu}
\date{7 September 2008}

\keywords{monomial ideal, Cohen--Macaulay, simplicial complex, local
cohomology, distraction, Alexander duality}

\maketitle

\section*{Introduction}

Scattered over the past few years have been several occurrences of
simplicial complexes whose topological behavior characterize the
Cohen--Macaulay property for quotients of polynomial rings by
arbitrary (not necessarily squarefree) monomial ideals.  It is unclear
whether researchers thinking about this topic have, to this point,
been aware of the full spectrum of related developments.  Therefore,
the purpose of this survey is to gather the developments into one
location, with self-contained proofs, including direct combinatorial
topological connections between them.

Four families of simplicial complexes are defined in reverse
chronological order: via distraction, the \v Cech complex, Alexander
duality, and then the Koszul complex.  Each comes with historical
remarks and context, including forays into Stanley decompositions,
standard pairs, $A$-hypergeometric systems, cellular resolutions, the
\v Cech hull, polarization, local duality, and duality for
$\ZZ^n$-graded resolutions.

Results or definitions appearing here for the first time include the
general categorical definition of cellular complex in
Definition~\ref{d:cellular}, as well as the statements and proofs of
Lemmas~\ref{l:M} and~\ref{l:fg}, though these are very easy.  The
characterization of exponent simplicial complexes in
Corollary~\ref{c:cover} might be considered new; certainly the
consequent connection in Theorem~\ref{t:exp} to the \v Cech simplicial
complexes is new, as is the duality in Theorem~\ref{t:dual} between
these and the dual \v Cech simplicial complexes.  The geometric
connection between distraction and local cohomology in
Theorem~\ref{t:V}, and its consequences in Section~\ref{s:geom},
generalize and refine results from \cite{BM}, in addition to providing
commutative proofs.  Finally, the connection between \v Cech and
Koszul simplicial complexes in Lemma~\ref{l:inclusion} and
Corollary~\ref{c:c-k}, as well as its general $\ZZ^n$-graded version
in Theorem~\ref{t:c-k}, appear~to~be~new.

\section{Monomial ideals and simplicial complexes}

Throughout this article, $I$ is an ideal in the polynomial ring $S =
\kk[x_1,\ldots,x_n]$.  The coefficient field~$\kk$ is assumed
arbitrary unless explicitly stated otherwise.  Each monomial in~$S$
can be expressed uniquely as $\xx^\aa = x_1^{a_1} \cdots x_n^{a_n}$
for some vector $\aa = (a_1,\ldots,a_n) \in \NN^n$ of $n$ nonnegative
integers.  More generally, Laurent monomials in
$S[x_1^{-1},\ldots,x_n^{-1}]$ can be expressed as $\xx^\aa$ for
vectors $\aa \in \ZZ^n$ of arbitrary integers.

For any subset $F \subseteq \{1,\ldots,n\}$ of indices, let $\NN^F
\subseteq \ZZ^F$ denote the \emph{coordinate subspace} of vectors
whose entries are zero outside of those indexed by~$F$.  This means
that $\ZZ^F$ consists of the vectors whose \emph{support}, by
definition the set $\supp(\aa)$ of indices~$i$ for which $a_i$ is
nonzero, is contained in~$F$.  Identify the set~$F$ of indices with
its characteristic vector in~$\NN^n$, which is the vector with
support~$F$ and all nonzero entries equal to~$1$.  Thus every
squarefree monomial is $\xx^F = \prod_{i \in F} x_i$ for some index
set $F$.  The inverse of such a monomial is~$\xx^{-F}$, so
\mbox{$S[\xx^{-F}] = S[x_i^{-1} \mid i \in F]$} is a \emph{monomial
localization} of~$S$, which we denote by~$S_F$.  In addition, let
\[
  \mm^F = \<\xx_F\> = \<x_i \mid i \in F\>
\]
be the ideal generated by the variables indexed by~$F$.  Thus $S/\mm^F
= \kk[\xx_\oF]$ is a polynomial ring in variables $\xx_\oF$ indexed by
the complementary set~\mbox{$\oF = \{1,\ldots,n\} \minus F$}.

For any ideal $I \subseteq S$, define the \emph{zero set}
\[
  Z(I) = \{\alpha \in \kk^n \mid f(\alpha) = 0 \text{ for all } f \in
  I\}.
\]
If $I$ is generated by monomials, for example, then $Z(I)$ is a union
of coordinate subspaces of~$\kk^n$.  Conversely, if any union~$Z$ of
coordinate subspaces is given, then the ideal $I(Z)$ of polynomials
vanishing on~$Z$ is a monomial ideal.  (When $\kk$ is finite, this
statement must be interpreted in the appropriate scheme-theoretic way;
but there is no harm in simply assuming that $\kk$ is algebraically
closed---and hence infinite---for the duration of this article.)  In
fact, $I(Z)$ is a \emph{squarefree} monomial ideal: its minimal
monomial generators are all squarefree.  On the other hand, a union
$Z$ of coordinate subspaces also gives rise to a \emph{simplicial
complex}
\[
  \Delta(Z) = \big\{F \subseteq \{1,\ldots,n\} \mid \kk^F \subseteq Z\big\},
\]
meaning that $\Delta(Z)$ is a collection of subsets, called
\emph{faces}, that is closed under inclusion: $F \in \Delta(Z)$ and $G
\subseteq F$ implies $G \in \Delta(Z)$.  The assignments $Z \goesto
I(Z)$ and $Z \goesto \Delta(Z)$ induce a bijection, called the
\emph{Stanley--Reisner correspondence}, between squarefree monomial
ideals and simplicial complexes.  Given a simplicial complex~$\Delta$,
the corresponding monomial ideal $I_\Delta$ is its
\emph{Stanley--Reisner ideal} (also known as its \emph{face ideal});
given a squarefree monomial ideal, $\Delta_I$ is its
\emph{Stanley--Reisner~complex}.

Background on squarefree monomial ideals and simplicial complexes can
be found in Chapter~II of Stanley's book \cite{Sta96}, or in Chapter~5
of the book by Bruns and Herzog \cite{BH93}, or in Chapter~1 of the
book by Miller and Sturmfels \cite{cca}.  The starting point for this
survey is the following.

\begin{defn}
A simplicial complex $\Delta$ is \emph{Cohen--Macaulay over\/~$\kk$}
if the Stanley--Reisner ring $\kk[\Delta]$ is Cohen--Macaulay.
\end{defn}

The Cohen--Macaulay condition on $\Delta$ is combinatorial, as
discovered by~\mbox{Reisner}.

\begin{thm}[Reisner's criterion, \cite{Rei76}]\label{t:reisner}
A simplicial complex $\Delta$ is Cohen--Macau\-lay over~$\kk$ if and
only if, for every face $F \in \Delta$, the subcomplex
\[
  \link F\Delta = \{G \in \Delta \mid G \cup F \in \Delta \hbox{ and }
  G \cap F = \nothing\},
\]
known as the \emph{link} of~$F$ inside $\Delta$, has reduced
cohomology
\[
  \HH^i(\link F\Delta;\kk) = 0\quad\hbox{for } i\neq\dim(\Delta)-|F|.
\]
\end{thm}
\begin{proof}
This follows from Theorem~\ref{t:hilb}, below; alternatively, see
\cite[Corollary~II.4.2]{Sta96}, \cite[Corollary~5.3.9]{BH93}, or
\cite[Section~13.5.2]{cca}.
\end{proof}

\begin{remark}\label{rk:top}
As combinatorial as Reisner's criterion may appear, in reality it is
topological: it is equivalent to the vanishing of all but the top
reduced cohomology of~$\Delta$ over~$\kk$ and the top \emph{local
cohomology} over~$\kk$ near every point $p$ in the geometric
realization $|\Delta|$; that is, the relative cohomology with
coefficients in~$\kk$ of the pair $(|\Delta|,|\Delta| \minus \{p\})$
\cite[Proposition~II.4.3]{Sta96}.
\end{remark}

\begin{remark}\label{r:char}
Note how the field $\kk$ enters into the definition: although the
presentation of the ring $\kk[\Delta]$ looks the same for every
field~$\kk$, the Cohen--Macaulay condition can depend on the
characteristic of~$\kk$ (but nothing else---all fields with the same
characteristic behave just like the prime field of that
characteristic).  This is reflected in Reisner's criterion, where the
simplicial complex is independent of~$\kk$, but its cohomology is
taken with coefficients~in~$\kk$.
\end{remark}

In recent years, Reisner's topological criterion for the
Cohen--Macaulay condition has been generalized to settings where the
monomial ideal~$I$ is no longer required to be squarefree.  The
essential framework is to define not one simplicial complex, but a
family of them indexed by subsets of~$\ZZ^n$, or ~$\kk^n$ in some
cases, and then connect the Cohen--Macaulayness of members in this
family to the homological algebra of~$I$, or to the geometry of
closely related ideals.  The purpose of this survey is to demonstrate
direct equivalences between various families of simplicial complexes
arising this way.

\section{Distracting arrangements and exponent complexes}\label{s:exp}

The first step is to produce a single family of simplicial complexes
whose Cohen--Macaulayness reflects that of a given monomial ideal~$I$,
which we fix from now on.

\begin{defn}\label{d:std}
The monomials in~$S$ outside of~$I$ are called \emph{standard
monomials}.  Write $\Lambda_I \subseteq \NN^n$ for the set of exponent
vectors on the standard monomials~of~$I$.
\end{defn}

\begin{lemma}\label{l:std}
The standard monomials form a $\kk$-vector space basis for $S/I$.\qed
\end{lemma}

The geometry of $\Lambda_I$ closely reflects the geometry of~$S/I$.
When $I$ is squarefree, for example, $\Lambda_I$ is the union of the
coordinate subspaces $\NN^F$ for the faces $F \in \Delta_I$ in the
Stanley--Reisner simplicial complex, while the zero set~$Z(I)$ is the
union of the coordinate subspaces $\kk^F$ for the faces $F \in
\Delta_I$.  In general, $\Lambda_I$ ``looks'' like the zero scheme
of~$I$, rather than the zero set~$Z(I)$.  Nonetheless, the transition
from standard monomials to subspace arrangement generalizes to
arbitrary monomial ideals, as long as $\kk$ has characteristic zero
(or characteristic bigger than the greatest exponent on any minimal
generator of~$I$).  For simplicity, we shall work with the field $\kk
= \CC$ of complex numbers.

\begin{defn}\label{d:V}
The \emph{distracting arrangement} $V_I$ of the monomial ideal~$I$ is
the Zariski closure $Z(I(\Lambda_I)) \subseteq \CC^n$ of the subset
$\Lambda_I \subseteq \NN^n \subseteq \CC^n$.
\end{defn}

Taking this Zariski closure has the same feel as the squarefree case.
For example, if a lattice point lies along a line containing
infinitely many points in~$\Lambda_I$, then the whole line is
contained in~$V_I$.  More generally, $V_I$ has the following
description.

\begin{lemma}\label{l:union}
$\Lambda_I$ can be expressed as a finite disjoint union of translates
$\aa + \NN^F\!$ of coordinate subspaces of\/~$\NN^n$.  For any such
decomposition, the distracting arrangement is the corresponding union
of translates $\aa + \CC^F\!$ \mbox{of coordinate
subspaces~of\/~$\CC^n$}.
\end{lemma}

The proof, and a fair bit of the rest of this survey, uses
$\ZZ^n$-gradings; see \cite[Chapter~8]{cca}.  In particular, $S =
\bigoplus_{\bb\in\ZZ^n} S_\bb$, where $S_\bb$ is the $1$-dimensional
vector space spanned by the monomial~$\xx^\bb$ if $\bb \in \NN^n$, and
$S_\bb = 0$ if $\bb \in \ZZ^n \minus \NN^n$.

\begin{proof}
Every $\ZZ^n$-graded $S$-module has a monomial associated prime
\cite[Proposition~8.11]{cca}, and hence a $\ZZ^n$-graded submodule
that is a $\ZZ^n$-graded translate $\xx^\aa\cdot S/\mm^F$ of a
quotient of~$S$ modulo a monomial prime.  By Noetherian induction,
$S/I$ possesses a filtration by $\ZZ^n$-graded submodules whose
associated graded pieces have this form.  The first sentence is thus a
consequence of Lemma~\ref{l:std}, which implies that the monomials
with exponent vectors in $\aa + \NN^F$ form a basis for $\xx^\aa\cdot
S/\mm^F$.  The second sentence follows from the first, since Zariski
closure commutes with finite~unions.\hspace{-2ex}%
\end{proof}

Unions as in Lemma~\ref{l:union} are called \emph{Stanley
decompositions} of~$S/I$.  They are relatively easy to come by,
although finding ones with good properties is difficult: Stanley's
conjecture \cite{Sta82}, which posits that one exists in which the
minimum size of any appearing~$F$ is at least the depth of~$S/I$, is
still open; see
\comment{}{\sf [ref], perhaps in this volume}\comment{} for recent
work on this topic and more \mbox{background}.

Locally near every point in~$\CC^n$, the distracting arrangement~$V_I$
looks like a coordinate subspace arrangement.  Here is a more precise
statement; see later in this section for history and context.

\begin{defn}\label{d:exp}
For $\alpha \in \CC^n$, let $Z_\alpha \subseteq V_I$ be the union of
the components of~$V_I$ passing through~$\alpha$.  Then $-\alpha +
Z_\alpha$ is a coordinate subspace arrangement in~$\CC^n$, and
$I(-\alpha + Z_\alpha)$ is the Stanley--Reisner ideal of the
\emph{exponent simplicial complex}~$\ED_\alpha$.
\end{defn}

\begin{remark}
Given a fixed set of monomial generators, the ideal $I \subseteq S$
can have different homological properties depending on the
characteristic of~$\kk$; see Remark~\ref{r:char}.  The simplicial
complexes~$\ED_\alpha$ for $\alpha \in \CC^n$, on the other hand, by
definition do not make reference to an arbitrary field~$\kk$.  It may
therefore come as a surprise that the topological invariants
of~$\ED_\alpha$ for $\alpha \in \ZZ^n$ control the homological
properties of~$S/I$ for any field~$\kk$; see Theorems~\ref{t:cech},
\ref{t:exp}, and~\ref{t:V}, as well as Corollary~\ref{c:V}.  But this
is analogous to the squarefree case, where the Stanley--Reisner
complex is defined without reference to~$\kk$, whereas its homological
invariants depend on~$\kk$; it justifies the restriction to $\kk =
\CC$ in Definition~\ref{d:V}.
\end{remark}


Although there is no known canonical way to choose a Stanley
decomposition, a measure of uniqueness can be imposed by dispensing
with the disjointness hypothesis.  As a matter of terminology, a
\emph{cover} of any set~$\Lambda$ is a family of subsets whose union
is~$\Lambda$.  One cover \emph{dominates} another if every subset from
the second cover is contained in a subset from the first.

\begin{prop}\label{p:cover}
$\Lambda_I$ has a unique cover by translates $\aa + \NN^F\!$ of
coordinate subspaces of\/~$\NN^n$ dominating every other such cover.
Each \emph{standard pair} \mbox{$\aa + \NN^F\!$} in this
\emph{standard pair decomposition} of~$\Lambda_I$ is the intersection
with $\Lambda_I$ of an irreducible component $\aa + \CC^F$ of the
distracting arrangement~$V_I$, and $\aa$ has
support~disjoint~from~$F$.
\end{prop}
\begin{proof}
Every component $\aa + \CC^F$ of~$V_I$ is the closure of a translated
orthant $\aa + \NN^F$ for some $\aa \in \NN^n$, by
Lemma~\ref{l:union}.  Observe that $(\aa + \CC^F) \cap \NN^n = \aa_F +
\NN^F$, where $\aa_F$ is obtained from~$\aa$ by setting to~$0$ the
coordinates of~$\aa$ indexed by~$\oF$.  Once $\Lambda_I$ contains $\aa
+ \NN^F$, it contains $\aa_F + \NN^F$, because the standard monomials
are closed under divisibility: every divisor of every standard
monomial is standard.
\end{proof}

Standard pair decompositions were introduced by Sturmfels, Trung, and
Vogel \cite[Section~3]{STV}, for the reason that the number of times
$F$ appears equals the multiplicity of the prime ideal $\mm^\oF$ in
the primary decomposition of~$I$.  Later, standard pairs became the
device through which the exponent simplicial complexes $\ED_\alpha$
were connected to the Cohen--Macaulay property.

\begin{remark}\label{rk:BM}
The earliest context in which exponent simplicial complexes seem to
have appeared explicitly, in the form of Definition~\ref{d:exp}, was
in the noncommutative algebra of multivariate hyper\-geometric systems
\cite{SST}.  The question of characterizing the Cohen--Macaulay
condition had arisen naturally.  The combinatorial commutative context
there was more general than the one here (it was for affine semigroup
rings), but it led back to distractions, nonetheless, as follows.

Associated to an affine semigroup is a family, parametrized
by~$\CC^d$, of systems of linear partial differential equations with
polynomial coefficients \cite{GGZ87,GKZ89}.  Based on observations by
Gelfand, Kapranov, and Zelevinsky, along with subsequent work by
Adolphson \cite{Ado94}, Sturmfels conjectured that the systems in this
family all have the same number of linearly independent solutions
precisely when the corresponding affine semigroup ring over the
complex numbers is Cohen--Macaulay.

This conjecture was proved by Matusevich, Miller, and Walther
\cite{MMW} by introducing functorial noncommutative Koszul-like
homological methods for $\ZZ^d$-graded modules over semigroup rings:
the rank of the solution space jumps whenever the parameter in~$\CC^d$
lies in a certain subspace arrangement. The arrangement is a
quasidegree set as in Corollary~\ref{c:exp}, coming from local
cohomology just as distractions come from monomial quotients.  In
particular, the arrangement is nonempty---so there exists at least one
rank jumping parameter in~$\CC^d$---precisely when there is nonzero
local cohomology to indicate the failure of Cohen--Macaulayness.

Berkesch and Matusevich \cite{BM} used an extension \cite{DMM} of
these homological methods to more general multigraded modules over
polynomial rings to see what the Cohen--Macaulay characterization by
rank jumps says about quotients by monomial ideals.  Their conclusions
were based on the fact that ranks of systems of differential equations
arising from monomial ideals are controlled by the geometry of their
distractions, from which Berkesch and Matusevich defined the exponent
simplicial complexes using standard pairs.  Sections~\ref{s:localCoh}
and~\ref{s:geom} contain details about the results in \cite{BM}; for
now, the reader is encouraged to consult that article for colorful
illustrations of distracting arrangements and exponent
simplicial~\mbox{complexes}.
\end{remark}

Proposition~\ref{p:cover} results in an alternate characterization of
exponent complexes that will be useful in the comparisons with other
simplicial complexes to~come.

\begin{cor}\label{c:cover}
A subset $F \subseteq \{1,\ldots,n\}$ is a maximal face (that is, a
\emph{facet}) of the exponent complex $\ED_\alpha$ of~$I$ if and only
if the set $(\alpha + \CC^F) \cap \Lambda_I$ of lattice points is
dense in $\alpha + \CC^F$ and $F$ is maximal with this property.
\end{cor}
\begin{proof}
$F$ is a facet of~$\Delta_\alpha$ if and only if $\alpha + \CC^F$ is
an irreducible component of the distracting arrangement~$V_I$, by
Definition~\ref{d:exp}.  If $\alpha + \CC^F$ is an irreducible
component of~$V_I$, then $(\alpha + \CC^F) \cap \Lambda_I$ is a dense
subset of it by Proposition~\ref{p:cover}, and $F$ is maximal with
this property because otherwise $\alpha + \CC^F$ would be strictly
contained in a com\-ponent $\alpha + \CC^{F'}$ for some other
face~$F'$.  On the other hand, if $(\alpha + \CC^F) \cap \Lambda_I$ is
dense in~$\alpha + \CC^F$, then $\alpha + \CC^F \subseteq V_I$ by
Definition~\ref{d:V}; maximality of~$F$ guarantees that $\alpha +
\CC^F$ is a component, given that each component of~$V_I$ is parallel
to a coord\-inate subspace of~$\CC^n$, which follows from
Proposition~\ref{p:cover},~or~even~Lemma~\ref{l:union}.
\end{proof}

\section{\v Cech simplicial complexes}\label{s:cech}

One of the many characterizations of the Cohen--Macaulay condition
\cite[Theorem~13.37]{cca} for a module~$M$ is by vanishing of all but
the top---that is, cohomological degree~$\dim(M)$---local
cohomology~$H^i_\mm(M)$ supported on the maximal ideal~$\mm$ of~$S$.
In particular, Reisner's criterion (Theorem~\ref{t:reisner}) is a
consequence of the following famous combinatorial formula for local
cohomology.

\begin{thm}[Hochster's formula]\label{t:hilb}
The $\ZZ^n$-graded Hilbert series of the local cohomology with maximal
support of a Stanley--Reisner ring satisfies
\[
  H(H^i_\mm(S/I_\Delta);\xx) = \sum_{F\in\Delta} \dim_\kk
  \HH^{i-|F|-1}(\link F\Delta;\kk) \prod_{j\in F}\frac{x_j^{-1}}{1 -
  x_j^{-1}}.
\]
\end{thm}
\begin{proof}
This can be deduced from Theorem~\ref{t:cech}, below; alternatively,
see \cite[Theorem~II.4.1]{Sta96}, \cite[Theorem~5.3.8]{BH93}, or
\cite[Theorem~13.13]{cca}.
\end{proof}

The proof of this statement is remarkably straightforward: the \v Cech
complex of $S/I_\Delta$ is a complex of $\ZZ^n$-graded $S$-modules,
whose graded piece in degree~$\aa$ is the cochain complex of the
appropriate link in~$\Delta$, with a homological shift.  This argument
is a quintessential example of the ``cellular'' technique that
pervades a significant portion of combinatorial commutative algebra in
recent years.  Quite generally, this technique begins with a cell
complex $X$, which could be a CW complex, a simplicial complex, a
polyhedral complex, or any other desired type of cell complex.  This
topological cell complex~$X$ has an algebraic chain complex
\[
\mathcal{C}_\spot X:\quad
0\ \from\hspace{-1ex}\bigoplus_{\text{vertices}\ v\in X}\ZZ_v 
 \ \from\hspace{-1ex}\bigoplus_{\text{edges}\ e\in X}\ZZ_e
 \ \from\cdots
 \ \from\hspace{-1ex}\bigoplus_{i\text{-faces}\ \sigma\in X}\ZZ_\sigma
 \ \from\cdots
\]
over the integers, where $\ZZ_\tau \from \ZZ_\sigma$ is multiplication
by some integer $\mathrm{coeff}(\sigma,\tau)$.

\begin{defn}\label{d:cellular}
Fix a cell complex~$X$, an abelian category, an object $J_\sigma$ for
each face $\sigma \!\in\! X$, and a natural map $J_\tau \!\from\!
J_\sigma$ for each inclusion $\tau \subseteq \sigma$
of\/~faces.~~The~\mbox{complex}
\[
J_\spot:\quad
0\ \from\hspace{-1ex}\bigoplus_{\text{vertices}\ v\in X}J_v 
 \ \from\hspace{-1ex}\bigoplus_{\text{edges}\ e\in X}J_e
 \ \from\cdots
 \ \from\hspace{-1ex}\bigoplus_{i\text{-faces}\ \sigma\in X}J_\sigma
 \ \from\cdots
\]
is \emph{cellular} and \emph{supported on}~$X$ if the component
$J_\tau \from J_\sigma$ is $\mathrm{coeff}(\sigma,\tau)$ times the
given natural map.  A similar construction yields \emph{cocellular}
complexes, starting from the algebraic cochain complex of~$X$ and the
natural maps are $J_\tau \to J_\sigma$ for~\mbox{$\tau \subseteq
\sigma$}.
\end{defn}

\begin{remark}
There are more variations than those mentioned in
Definition~\ref{d:cellular}.  For example, we could start from the
\emph{reduced} (or \emph{augmented}) algebraic chain or cochain
complex of a simplicial complex, which includes another
(co)homological degree for the empty face of~$X$.  We could also start
from the algebraic relative (co)chain complex of a pair of cell
complexes.  Generally speaking, any complex in an abelian category
constructed from underlying topological data on a cell complex, or
pair of cell complexes, is referred to as a \emph{cellular~complex} of
objects in that~category.
\end{remark}

\begin{example}\label{ex:cellular}
The prototypical example is the algebraic cochain complex arising from
\v Cech cohomology theory in topology.  In this case, the cell
complex~$X$ is a simplex, and $J_\sigma$ is a space of functions on an
intersection $\bigcap_{j \in \sigma} U_j$ of open sets in a given
cover $\{U_1,\ldots,U_n\}$ of a fixed topological space (usually
unrelated to~$X$).

In its $\ZZ^n$-graded avatar, the \emph{\v Cech complex}
$\vC^\spot(\xx)$ is a cocellular complex of monomial localizations
of~$S$ supported on a simplex~$X$.  More precisely, $X$ is the simplex
on~$\{1,\ldots,n\}$, and $J_F = S_F$ is the monomial localization
$S[\xx^{-F}]$, the natural map $S_F \to S_{F'}$ being the localization
homomorphism for $F \subseteq F'$.
\end{example}

\begin{remark}
Other cellular complexes have played roles in commutative algebra for
well over half a century, perhaps the first being the Koszul complex
\cite{koszul}, which was originally introduced for the purpose of Lie
algebra cohomology theory.  The initial use of cellular complexes in
identifiably combinatorial commutative algebra came with the Taylor
complex \cite{taylor}, which provides a free resolution of any
monomial ideal in a polynomial ring.  The application of \v Cech
complexes to local cohomology of monomial ideals was carried out in
the 1970's.  The earliest cellular complexes supported on polyhedra
that are not simplicial complexes arrived in the work of Ishida on
dualizing complexes and local cohomology with maximal support over
affine semigroup rings \cite{Ish80,Ish87}; in addition, Ishida's
objects~$J_\sigma$ were modules over rings other than polynomial
rings.  The notion of \emph{cellular resolution} in combinatorial
commutative algebra was codified formally for the first time in the
work of Bayer and Sturmfels, also with Peeva \cite{BS,BPS}, in the
context of lattice ideals and generic monomial ideals (the more
inclusive definition of generic monomial ideals in current use
\cite{MSY} was introduced later).

More general objects~$J_\sigma$, including graded injectives and
quotients by irreducible monomial ideals in polynomial and affine
semigroup rings, began to appear in \cite{Mil00} and \cite{MilIrr},
where topological duality for the support complexes was observed to
reflect algebraic (Matlis or local) duality algebraically.  All of the
cellular resolutions to that point had been supported on special
classes of cell complexes, notably simplicial and polyhedral
complexes.  With the work of Batzies and Welker \cite{BW}, more
flexible support complexes~$X$ began appearing, allowing the
application of homotopic and Morse-theoretic cancellative approaches,
which don't preserve polyhedrality.  Cellular methods generally bring
topological or combinatorial techniques to bear on algebraic problems,
but Fl\o ystad, with his \emph{enriched cohomology}, has begun to
realize that commutative algebra can provide additional homological
insight into the cell complex~$X$, rather than vice versa
\cite{floystad}.  Most recently, cellular techniques in combinatorial
commutative algebra have given way to non-combinatorial applications
where the objects~$J_\sigma$ do not belong to commutative algebra at
all: in \cite{multIdealSums}, where Definition~\ref{d:cellular} was
mentioned in the Introduction without details, they are sheaves from
scheme-theoretic algebraic geometry, or perhaps from complex analytic
geometry.

These examples are idiosyncratic choices of landmarks; there have been
many others between, demonstrating the strengths and limitations of
cellular \mbox{resolutions}.
\end{remark}

Generalizing Theorem~\ref{t:hilb} to the local cohomology of $S/I$ for
an arbitrary monomial ideal~$I$ involves simplicial complexes that can
be characterized as follows.

\begin{defn}\label{d:cech}
Given $\aa \in \ZZ^n$, write $\aa = \aa_+ - \aa_-$ as a difference of
vectors in~$\NN^n$.  The \emph{\v Cech simplicial complex} of~$I$ in
degree~$\aa$ is the set
\[
  \CD_\aa = \{F - \supp(\aa_-) \mid (S_F/IS_F)_\aa = \kk\}
\]
in bijection with those $F$ for which the localization $S/I \otimes_S
S_F$ is nonzero in degree~$\aa$.
\end{defn}

\begin{lemma}
$\CD_\aa$ is a simplicial complex.
\end{lemma}
\begin{proof}
For $F \subseteq F'$, the homomorphism from $S_F/IS_F$ to its image in
$S_{F'}/IS_{F'}$ is the quotient homomorphism of $S_F/IS_F$ modulo the
kernel of the localization that inverts~$\xx^{F'-F}$.  Therefore, once
$S_F$ is nonzero in degree~$\aa$, the homomorphism $(S_F/IS_F)_\aa \to
(S_{F'}/IS_{F'})_\aa$ is a quotient map.  The final observation to
make is that $S_F$ is nonzero in degree~$\aa$ precisely when $F$
contains~$\supp(\aa_-)$.
\end{proof}

\begin{thm}\label{t:cech}
The \v Cech simplicial complexes $\CD_\aa$ of~$I$ compute local
cohomology:
\[
  H^i_\mm(S/I)_\aa = \HH^{i-|\supp(\aa_-)|-1}(\CD_\aa;\kk).
\]
\end{thm}
\begin{proof}
The local cohomology $H^i_\mm(S/I)_\aa$ is the $\ZZ^n$-degree $\aa$
piece of the \v Cech complex $S/I \otimes_S \vC^\spot(\xx)$.  Viewing
the \v Cech complex as being cocellular as in
Example~\ref{ex:cellular}, and using that $S/I \otimes_S S_F =
S_F/IS_F$, the degree~$\aa$ piece of $S/I \otimes_S \vC^\spot(\xx)$ is
the reduced cochain complex of~$\CD_\aa$ by Definition~\ref{d:cech},
with the faces $F - \supp(\aa_-)$ of dimension~$i - |\supp(\aa_-)| -
1$ contributing to cohomological degree $i = |F|$.
\end{proof}

Theorem~\ref{t:cech} was recorded by Takayama
\cite[Theorem~1]{takayama}, as the first step in studying locally
Cohen--Macaulay monomial schemes---that is, whose local cohomology
modules $H^i_\mm(S/I)$ have finite length for $i < \dim(S/I)$.
Takayama also proved, using the simplicial formula, that
$H^i_\mm(S/I)$ is nonzero only in degrees $\aa \in \ZZ^n$ satisfying
$\aa \preceq \aa_I - \1$, where $\aa_I$ is the exponent on the least
common multiple of the minimal generators, $\1 = (1,\ldots,1)$, and
the relation $\preceq$ is coordinatewise dominance.  This obsesrvation
admits a simple proof in the general context of arbitrary finitely
generated $\ZZ^n$-graded modules by using $\ZZ^n$-graded injective
resolutions (see \cite[Chapter~11]{cca} for an introduction).  The
general statement is as follows.

\begin{lemma}\label{l:M}
If the Betti numbers $\beta_{i,\bb}(M)$ of a finitely generated
$\ZZ^n$-graded module~$M$ lie in degrees $\bb \preceq \aa$, then
$H^i_\mm(M)$ can only be nonzero in $\ZZ^n$-degrees $\preceq \aa -
\1$.
\end{lemma}
\begin{proof}
The proof consists of two~\mbox{observations}:
\begin{enumerate}[(i)]
\item
the Betti number $\beta_{i,\bb}(M)$ equals the Bass
number~$\mu_{n-i,\bb-\1}(M)$ at the maximal ideal, as can be seen by
using a Koszul complex to compute both; and
\item
the $\ZZ^n$-graded injective hull of~$\kk$ is nonzero only in degrees
from~$-\NN^n$.\qedhere
\end{enumerate}
\end{proof}

In fact, (i)~implies more, when $M$ is graded by~$\NN^n$, and not
just~$\ZZ^n$.

\begin{lemma}\label{l:fg}
If $M$ is $\NN^n$-graded and $H^i_\mm(M)$ is finitely generated (i.e.,
has finite length), then the $\ZZ^n$-graded translate
$H^i_\mm(M)(-\1)$ up by $(1,\ldots,1)$ is also~$\NN^n$-graded.
\end{lemma}
\begin{proof}
The \emph{\v Cech hull} \cite[Definition~2.7]{Mil00} (or see
\cite[Definition~13.32]{cca}) of any $\ZZ^n$-graded module~$N$ is the
$\ZZ^n$-graded module $\vC N$ whose degree~$\bb$ piece is $(\vC N)_\bb
= N_{\bb_+}$ for $\bb \in \ZZ^n$.  The \v Cech hull is an exact
functor on $\ZZ^n$-graded modules because this is measured
degree-by-degree.  It is easy to show that the \v Cech hull of any
indecomposable $\ZZ^n$-graded injective module~$E$ is either zero
or~$E$, the latter precisely when the $\ZZ^n$-graded degree piece
$E_\0$ is nonzero \cite[Lemma~4.25]{Mil00}.  It follows immediately
from~(i) that the minimal injective resolution of~$M(-\1)$---and hence
all of the local cohomology of~$M(-\1)$---is fixed by the \v Cech
hull.  The previous sentence is what the statement of the lemma should
really say, since any finitely generated module fixed by the \v Cech
hull must vanish in all degrees~\mbox{$\not\in \1 + \NN^n$}.
\end{proof}

\begin{remark}
Rahimi generalized Theorem~\ref{t:cech} to a combinatorial formula for
the local cohomology of $S/I$ supported at an arbitrary prime monomial
ideal \cite{rahimi}.  It remains open to find a combinatorial formula
for the local cohomology of $S/I$ supported at an arbitrary (may as
well be squarefree) monomial ideal.
\end{remark}

\section{Local cohomology vs.\ distraction}\label{s:localCoh}

Thus far we have seen two families of simplicial complexes: the
exponent (Definition~\ref{d:exp}) and \v Cech
(Definition~\ref{d:cech}) complexes.  They were introduced in
radically different contexts, the former being the geometry of
distractions, and the latter being the graded algebra of local
cohomology; but they are very closely related.

\begin{thm}\label{t:exp}
The exponent simplicial complex~$\ED_\aa$ of~$I$ in degree $\aa \in
\ZZ^n$ is a cone with apex $\supp(\aa_-)$, and the \v Cech simplicial
complex $\CD_\aa$ is the link of~the~apex:
\[
  \CD_\aa = \link{\supp(\aa_-)}{\ED_\aa}.
\]
More generally, for $\alpha \in \CC^n$, define \mbox{$\alpha' \in
\ZZ^n$} by setting to~$-1$ all coordinates of~$\alpha$ that lie
outside of\/~$\NN$.  Then $\ED_\alpha$ is a cone with apex~$\alpha'_-$
and \mbox{$\CD_{\alpha'} = \link{\alpha'_-}{\ED_\alpha}$}.
\end{thm}
\begin{proof}
First observe that $\CD_{\aa'} = \CD_\aa$ for $\aa \in \ZZ^n$, since
the variables $x_i$ for $i \in \supp(\aa_-) = \aa'_-$ are units on
$S_F/IS_F$ whenever $F \supseteq \supp(\aa_-)$.  Hence the second
sentence is indeed more general than the first.

\begin{lemma}\label{l:star}
$\ED_{\alpha'} = \ED_\alpha$ for all $\alpha \in \CC^n$.  More
generally, if~$\beta$ is obtained from $\alpha \in \CC^n$ by changing
some nonnegative integer entries of~$\alpha$ to lie outside
of\/~$\NN$, then $\ED_\beta = \Star
{\supp(\alpha-\beta)}{\ED_\alpha}$, where for any simplicial
complex~$\Delta$ and face~$F$,
\[
  \Star F\Delta = \{F' \in \Delta \mid F \cup F' \in \Delta\}.
\]
\end{lemma}
\begin{proof}
The first sentence follows from the second by expressing both as stars
in $\ED_{\alpha_0}$, where $\alpha_0$ is obtained from $\alpha$ by
setting to~$0$ all entries of~$\alpha$ lying outside~of~$\NN$.

For $F$ to be a facet of $\ED_\alpha$, it must by
Corollary~\ref{c:cover} contain the set of indices $i$ for which
$\alpha_i$ lies outside of~$\NN$, or else the affine subspace $\alpha
+ \CC^F$ fails to intersect~$\NN^n$, let alone being the closure of a
set of lattice points therein.  Therefore the facets of~$\ED_\beta$
are precisely the facets of~$\ED_\alpha$ that contain
$\supp(\alpha-\beta)$.
\end{proof}

The lemma reduces Theorem~\ref{t:exp} to the claim that $\CD_\aa =
\link{\aa_-}{\ED_\aa}$ when \mbox{$\aa' = \aa$}.  For this, it is
enough by Corollary~\ref{c:cover} to show, when $\aa_- \subseteq F$,
that $(S_F/IS_F)_\aa = \kk$ if and only if $(\aa + \CC^F) \cap
\Lambda_I$ is dense in $\aa + \CC^F$; but both of these, given $\aa_-
\subseteq F$, are equivalent to $\xx^{\aa_+}\CC[\xx_F]$ being a
subquotient of~$S/I$.
\end{proof}

\begin{remark}\label{rk:V}
Berkesch and Matusevich observed (cf.\ Remark~\ref{rk:BM}) that the
cohomology of the exponent simplicial complexes of~$I$ controls the
Cohen--Macaulay property.  Their precise statement along these lines
\cite[Theorem~3.14]{BM} is the final claim in Corollary~\ref{c:V},
below.  However, they did not explicitly mention that the
Cohen--Macaulay property for all exponent simplicial complexes is
equivalent to the Cohen--Macaulay property for the distracting
arrangement~$V_I$, because their definition of $\ED_\alpha$ was
combinatorial, in contrast to the equivalent geometric definition in
Section~\ref{s:exp}~here.  Nevertheless, the next main result here,
Theorem~\ref{t:V}, which provides a precise geometric criterion for
the failure of the Cohen--Macaulay condition at a point of~$V_I$, is
inspired by a similar---but weaker---result of Berkesch and
Matusevich; see \mbox{Corollary~\ref{c:exp}}.
\end{remark}

\begin{defn}
Given a $\ZZ^n$-graded module~$N$, its \emph{true degree set} is the
set
\[
  \tdeg(M) = \{\aa \in \ZZ^n \mid N_\aa \neq 0\}
\]
of $\ZZ^n$-graded degrees where $N$ is nonzero.  The \emph{quasidegree
set} $\qdeg(N)$ is the Zariski closure of~$\tdeg(N)$ inside
of~$\CC^n$.  Let $\qdeg_\RR(N) = \qdeg(N) \cap \RR^n$ be the
\emph{real quasidegree set} of~$N$.
\end{defn}

\begin{example}
The true degree set of $S/I$ is $\tdeg(S/I) = \Lambda_I$ from
Definition~\ref{d:std}.  The quasidegree set of $S/I$ is the
distracting arrangement~$V_I$ from Definition~\ref{d:V}.
\end{example}

\begin{remark}
Quasidegree sets and distracting arrangements in~$\CC^n$ are varieties
over $\CC$, so it makes sense to say they are Cohen--Macaulay.  Their
intersections with~$\RR^n$ are varieties over $\RR$, but are also
topological spaces.  Using Remark~\ref{rk:top}, it therefore makes
sense to say that quasidegree sets or distracting arrangements
over~$\RR$ are Cohen--Macaulay over a field that might not be~$\CC$.
That said, the following result is true---with the same proof---if
$\RR$ and~$\kk$ are replaced by~$\CC$ throughout.
\end{remark}

\begin{thm}\label{t:V}
If $\dim(S/I) = d$, the \emph{real distracting arrangement} $V_I(\RR)
= V_I \cap \RR^n$ is not Cohen--Macaulay over\/~$\kk$ near $\alpha \in
\RR^n \iff \alpha \in
\qdeg_\RR\!\big(\hspace{-.062ex}H^i_\mm(S/I)\hspace{-.062ex}\big)$
for~some~\mbox{$i < d$}.
\end{thm}
\begin{proof}
The arrangement $V_I(\RR)$ fails to be Cohen--Macaulay over~$\kk$
near~$\alpha \in \RR^n$ if and only if $\ED_\alpha$ fails to be
Cohen--Macaulay over~$\kk$.
Noting that $\link F{\ED_\alpha}$ remains a cone unless $F$ contains
$\alpha'_-$ (Theorem~\ref{t:exp}), this occurs precisely when there
exist a face $F \in \ED_\alpha$ with $F \supseteq \alpha'_-$ and a
cohomological degree $i < d$~such~that
\[
  \HH^{i-|F|-1}(\link F{\ED_\alpha};\kk) \neq 0
\]
by Reisner's criterion (Theorem~\ref{t:reisner}).  On the other hand,
if $\aa = \alpha_F$ is the result of setting to~$-1$ all coordinates
of~$\alpha$ indexed by~$F$, then $\link F{\ED_\alpha} = \link
F{\ED_\aa}$ by Lemma~\ref{l:star}.  But $\link F{\ED_\aa} = \CD_\aa$
by Theorem~\ref{t:exp}.  Furthermore, $\CD_\aa = \CD_\bb$ for all $\bb
\in \aa - \NN^F$ by the first sentence of the proof of
Theorem~\ref{t:exp}.  Stringing these statements together, and
applying Theorem~\ref{t:cech}, we find that $V_I(\RR)$ fails to be
Cohen--Macaulay over~$\kk$ at~$\alpha$ if and only if
\[
\tag{$*$}
  \text{there exist }F\supseteq\alpha'_-\text{ and }i < d\text{ with }
  \bb \in \tdeg(H^i_\mm(S/I)) \text{ for all }\bb \in\alpha_F - \NN^F.
\]
This condition certainly implies that $\alpha \in \qdeg_\RR
H^i_\mm(S/I)$.

For the converse, assume $i < d$ with $\alpha \in \qdeg_\RR
(H^i_\mm(S/I))$.  Then $\alpha + \RR^F$ contains a dense subset of
points in $\tdeg(H^i_\mm(S/I))$ for some~$F$ and $i < d$.  This set
$F$ is forced to contain the set of indices where $\alpha$ is not an
integer, or else $\alpha + \RR^F$ wouldn't contain any lattice points
at all.  Replacing $\alpha$ by some point in $\alpha + \RR^F$ with
$\alpha_i = -1$ for $i \in F$ (this changes neither $\alpha_F$
nor~$\alpha'$), we may as well assume $\alpha \in \ZZ^n$.
Lemma~\ref{l:M} implies that $\tdeg(H^i_\mm(S/I)) \cap (\alpha -
\NN^F)$ contains a lattice point~$\bb$, since the Zariski closure of
the complementary set $\tdeg(H^i_\mm(S/I)) \cap (\alpha - (\ZZ^F
\minus \NN^F))$ of lattice points in $\alpha + \ZZ^F$ has dimension
less than~$|F|$.  Theorems~\ref{t:cech} and~\ref{t:exp} imply that
$\HH^{i-|\bb'_-|-1}(\link {\bb'_-}{\ED_\bb};\kk) \neq 0$.  Given this
nonvanishing, we might as well replace $F$ with $\bb'_-$, which
contains $\alpha'_-$ by construction, proving~($*$).
\end{proof}

\begin{remark}
The analysis in Theorem~\ref{t:V} and its proof comes down to the fact
that if some exponent simplicial complex $\ED_\alpha$ fails to be
Cohen--Macaulay, then some link in $\ED_\alpha$ has nonzero reduced
cohomology in cohomological degree less than its dimension.  But the
\v Cech simplicial complex~$\CD_\alpha$, which is indeed a link in
$\ED_\alpha$, might not itself have nonzero non-top reduced
cohomology.  Thus, even if you're standing at a lattice point $\alpha$
whose simplicial complex fails Cohen--Macaulayness, you might not get
nonvanishing local cohomology in that particular degree; instead, you
might have travel along a component of $\qdeg(H^i_\mm(S/I))$, thereby
taking a link, to get to a lattice point whose \v Cech simplicial
complex does indeed have nonzero non-top local cohomology.  That's the
point, really: the difference between the distracting
arrangement~$V_I$ and $\qdeg(H^i_\mm(S/I))$ is the difference between
\emph{possessing} a face whose link witnesses failure of
Cohen--Macaulayness and actually \emph{choosing}~one.
\end{remark}

\section{Geometric Cohen--Macaulay criteria}\label{s:geom}

The equivalence between exponent and \v Cech simplicial complexes in
Theorem~\ref{t:exp}, and the geometric relation between the
distracting arrangement and the degrees of nonvanishing local
cohomology in Theorem~\ref{t:V}, give rise to a number of
Cohen--Macaulay criteria.  Here is the first, already discussed in
Remark~\ref{rk:V}.

\begin{cor}\label{c:V}
$S/I$ is Cohen--Macaulay
$\iff$ the real distracting arrangement~$V_I(\RR)$ is Cohen--Macaulay
over\/~$\kk$.
In fact, $S/I$ is Cohen--Macaulay of dimension~$d$
$\iff$ the exponent complexes $\ED_\alpha$ for $\alpha \in \CC^n$ are
Cohen--Macaulay over\/~$\kk$ of dimension~$d - 1$.
\end{cor}
\begin{proof}
By the well-known local cohomology criterion for Cohen--Macaulayness
(see \cite[Theorem~13.37.9]{cca}, for example), Theorem~\ref{t:V}, and
Definition~\ref{d:exp},
\begin{align*}
S/I \text{ is Cohen--Macaulay}
  &\iff H^i_\mm(S/I) = 0 \text{ for all } i < d
\\&\iff \qdeg_\RR(H^i_\mm(S/I))
        \text{ is empty for all } i < d
\\&\iff V_I(\RR) \text{ is Cohen--Macaulay over }\kk \text{ at every
        point } \alpha\!\in\!\RR^n
\\&\iff \ED_\alpha\!\text{ is Cohen--Macaulay over }\kk \text{ at
        every point } \alpha\in\RR^n.
\end{align*}
The result now follows from Lemma~\ref{l:star}, which implies that
$\{\ED_\alpha \mid \alpha \in \CC^n\} = \{\ED_\alpha \mid \alpha \in
\RR^n\} = \{\ED_\aa \mid \aa \in \ZZ^n\}$; that is, no new exponent
simplicial complexes arise by considering $\alpha \in \CC^n$ instead
of $\alpha \in \RR^n$.
\end{proof}

This connection between the geometry of distracting arrangements and
homological properties of monomial ideals can also be seen to arise
from elementary algebraic manipulation of monomials.

\begin{defn}
For a variable $t$ and $a \in \NN$, the \emph{distraction} of the pure
power $t^a$~is
\[
  \wt{t^a} = t(t-1)(t-2)\cdots(t-a+1),
\]
a polynomial of degree~$a$; the distraction of $1 = t^0$ is defined to
equal~$1$.  For $\aa \in \NN^n$, the \emph{distraction} of the
monomial $\xx^\aa \in S$ is the product
\[
  \wt{\xx^\aa} = \prod_{i=1}^n \wt{x_i^{a_i}}
\]
of the distractions of its pure-power factors~$x_i^{a_i}$.  For any
monomial ideal $I \subseteq S$, the \emph{distraction} of~$I$ is the
ideal~$\wt I$ generated by the distractions of the monomials in~$I$.
\end{defn}

\begin{example}
The distraction of $I = \<x^3yz^2, y^2z^4\>$ is
\[
  \wt I = \<x(x-1)(x-2)yz(z-1),\ y(y-1)z(z-1)(z-2)(z-3)\>.
\]
\end{example}

\begin{lemma}
$\wt I$ is generated by the distractions of the minimal generators
of~$I$.\qed
\end{lemma}

The definitions allow for an algebraic version of Corollary~\ref{c:V}.

\begin{cor}\label{c:V'}
If\/ $\kk = \CC$, then $S/I$ is Cohen--Macaulay $\iff$ the distraction
$S/\tilde I$ is.
\end{cor}
\begin{proof}
$S/\tilde I$ is the complex coordinate ring of the distracting
arrangement~$V_I$.
\end{proof}

Why should the algebraic formulation of distraction lead one to think
that it detects Cohen--Macaulayness?  It ought to because regular
sequences connect a single object, the polarization, to both $S/I$ and
$S/\tilde I$.  This connection has been observed before; it was noted
explicitly in the Introduction to \cite{BCR}, for instance.

\begin{defn}
For a variable $t$ and $a \in \NN$, the \emph{polarization} of the
pure power~$t^a$~is
\(
  t'_a = t_0 t_1 t_2 \cdots t_{a-1},
\)
a squarefree monomial in $a$ new variables $t' = t_0,\ldots,t_{a-1}$.
\pagebreak
For $\aa \in \NN^n$, the \emph{polarization} of the monomial $\xx^\aa
\in S$ is the product
\[
  {\xx'}_{\!\!\aa} = \prod_{i=1}^n (x'_i)_{a_i} \in S'
\]
of the polarizations of its pure-power factors~$x_i^{a_i}$, where $S'$
is a polynomial ring in variables $x_{ij}$ for $i = 1,\ldots,n$ and $j
= 0,\ldots,a_i-1$.  The \emph{polarization} $I' \subseteq S'$ of a
monomial ideal $I \subseteq S$ is generated by the polarizations of
its minimal generators.
\end{defn}

\begin{example}
The polarization of $I = \<x^3yz^2, y^2z^4\>$ is
\[
  I' = \<x_0x_1x_2 y_0 z_0z_1,\ y_0y_1 z_0z_1z_2z_3\>.
\]
\end{example}

It is well-known that the quotient by a monomial ideal~$I$ is
Cohen--Macaulay if and only if the quotient by its polarization~$I'$
is; see \cite[Exercise~3.15]{cca} or
\cite[Proposition~1.3.4]{kummini08}.  The idea of the proof is that
the elements $x_{ij} - x_{i0}$ form a regular sequence in~$S'/I'$, and
the quotient modulo this regular sequence is~$S/I$, if one sets
$x_{i0} = x_i$ for all~$i$.  This argument in fact shows that a
minimal free resolution of $S'/I'$ over~$S'$ descends to a minimal
free resolution of $S/I$ over~$S$, so most of the homological
invariants---Betti numbers, and so on---descend, as well.  A similar
procedure relates the polarization to the distraction, so
Corollary~\ref{c:V'} is not surprising from this perspective, although
it is harder than one might imagine to prove it directly via regular
sequences.  Indeed, $S/\tilde I$ is not a local ring, so the
implication ``$V_I$ is Cohen--Macaulay $\implies S/I'$ is
Cohen--Macaulay'' isn't obvious.

The geometric characterization of the Cohen--Macaulay condition via
distracting arrangements in Corollary~\ref{c:V} gives rise to another
geometric characterization, also observed by Berkesch and Matusevich
\cite[Corollary~1.4]{BM}: when $V_I$ is sliced by all translates of a
sufficiently general linear subspace of dimension $n-\dim(S/I)$, the
number of points in the intersection---counted with multiplicity---is
constant precisely when $S/I$ is Cohen--Macaulay.  The following
generalizes their result, which had additional hypotheses (such as
integrality) on the projection~$A$.

\begin{prop}\label{p:deg}
Let $\kk = \CC$.  If a linear map $A: \CC^n \to \CC^d$ takes every
coordinate subspace in the zero set $Z(I)$ isomorphically to~$\CC^d$,
then every fiber of the projection \mbox{$V_I \to \CC^d$} has Krull
dimension~$0$.  $S/I$ is Cohen--Macaulay of dimension~$d$ exactly when
every fiber has multiplicity~$\deg(I)$ equal to the standard
$\ZZ$-graded~degree~of~$I$.
\end{prop}
\begin{proof}
$S/I$ is Cohen--Macaulay if and only if $V_I$ is, by
Corollary~\ref{c:V} along with Corollary~\ref{c:V'} and its proof.
The projection~$A$ makes $V_I$ into a family over the smooth
base~$\CC^d$.  Since $V_I$ is a finite union of subspaces parallel to
coordinate spaces in~$Z(I)$, this family is finite over~$\CC^d$.
Under these conditions, the Cohen--Macaulay condition on~$V_I$ is
equivalent to the flatness of this family (see
\cite[Proposition~2.2.11]{BH93} or \cite[Theorem~13.37.5]{cca}), which
in turn is equivalent to the constancy of the multiplicity of the
fibers \cite[Proposition~III.9.2(e) and Exercise~II.5.8]{Har}.  The
constant multiplicity must be $\deg(I)$ because this is the number of
dimension~$d$ irreducible components of~$V_I$.
\end{proof}

The statement from which Berkesch and Matusevich derive most of their
consequences, such as the second half of Corollary~\ref{c:V} and their
version of Proposition~\ref{p:deg}, is a coarsening
\cite[Theorem~1.3]{BM} of Theorem~\ref{t:V} to a $\ZZ^d$-graded
setting.  The situation is that of Proposition~\ref{p:deg}, except
that $A$ is required to be represented by a $d \times n$ integer
matrix, also called~$A$, whose columns span~$\ZZ^d$ and whose kernel
contains no nonzero positive vectors.  Such a matrix~$A$ induces a
$\ZZ^d$-grading on~$S$, in which each monomial $\xx^\aa$ has degree
$A\cdot\aa \in \ZZ^d$.  Any module $N$ that is $\ZZ^d$-graded has
\emph{true degrees} $\tdeg_d(N) = \{\uu \in \ZZ^d \mid N_\uu \neq 0\}$
and \emph{quasidegrees} $\qdeg_d(N) = $ Zariski closure of
$\tdeg_d(N)$ in~$\CC^d$.

\begin{lemma}\label{l:qdeg}
Fix a $\ZZ^n$-graded finitely generated $S$-module $M$ of dimension at
most~$d$, or an artinian $\ZZ^n$-graded module whose $\ZZ^n$-graded
Matlis dual has dimension at most~$d$.  Then $\qdeg_d(M) =
A\big(\!\qdeg(M)\big)$ is the image in~$\CC^d$ of
the~\mbox{$\CC^n$-quasidegrees}.
\end{lemma}
\begin{proof}
By Matlis duality, assume $M$ is finitely generated.  The set
$\tdeg_n(M)$ of true degrees of~$M$ admits a decomposition as a finite
union of positive integer orthants having the form $\alpha + \NN^F$,
where $\alpha \in \ZZ^n$ and $|F| \leq d$.  Any such decomposition
yields a decomposition of $\tdeg_d(M)$ as a finite union of translated
semigroups $\beta + \NN A_F$, where $A_F$ is the set of columns of~$A$
indexed by~$F$.  The result holds because Zariski closure commutes
with the application of~$A$ to such sets.
\end{proof}

Here, finally, is a precise formulation of \cite[Theorem~1.3]{BM}.

\begin{cor}\label{c:exp}
Fix a $\ZZ^d$-grading by~$A$.  In the situation of
Proposition~\ref{p:deg}, the fiber over $\gamma \in \CC^d$ has
multiplicity $> \deg(I)$ exactly when $\gamma \in
\bigcup_{i=0}^{d-1}\qdeg_d(H^i_\mm(S/I))$.
\end{cor}
\begin{proof}
The multiplicity of the fiber $V_I(\gamma)$ over $\gamma \in \CC^d$ is
$> \deg(I)$ if and only if $V_I$ fails to be flat near $\gamma$ as a
family over~$\CC^d$.  This occurs precisely when the stalk of the
coordinate ring $\CC[V_I] = S/\tilde I$ near~$\gamma$ fails to be
Cohen--Macaulay over the local ring of $\gamma \in \CC^d$.  This
condition is equivalent to the failure of the Cohen--Macaulay
condition for $V_I$ locally at some point $\alpha \in \CC^n$ mapping
to~$\gamma$ under~$A$.  But we have already seen, in
Theorem~\ref{t:V}, that the set of points in~$V_I$ where the
Cohen--Macaulay condition fails is the union, over $i = 0,\ldots,d-1$,
of the quasidegree sets $\qdeg(H^i_\mm(S/I))$ in~$\CC^n$.  The result
now follows from Lemma~\ref{l:qdeg}.
\end{proof}

\section{Dual \v Cech simplicial complexes}\label{s:dual}

What would a survey about topological criteria for Cohen--Macaulayness
be without some mention of Alexander duality?  Well, here it is, in a
third family of simplicial complexes that has appeared in the
literature and whose cohomology detects the Cohen--Macaulay property
for an arbitrary monomial ideal~$I$.

\begin{defn}\label{d:dual}
Write $\ol\Lambda_I = \NN^n \minus \Lambda_I$ for the set of exponent
vectors on the monomials in~$I$.  Given $\bb \in \ZZ^n$, write $\bb =
\bb_+ - \bb_-$ as a difference of vectors in~$\NN^n$.  The \emph{dual
\v Cech simplicial complex} for $\aa \in \ZZ^n$ is
\[
  \CD^\aa = \big\{F \subseteq \{1,\ldots,n\} \mid F = \supp(\bb_-)
  \text{ for some } \bb \in \aa + \ol\Lambda_I \big\}.
\]
\end{defn}

As suggested by the nomenclature, Definition~\ref{d:dual} is connected
to \v Cech simplicial complexes by Alexander duality, as we shall see
later in this section.

\begin{defn}\label{d:alexdual}
The \emph{Alexander dual}\/ $\Delta^*$ of a simplicial
complex~$\Delta$ is characterized by the following property: $F
\not\in \Delta$ if and only if $\oF \in \Delta^*$.
\end{defn}

Alexander duality for simplicial complexes entered into the earliest
of Hochster's investigations into monomial commutative algebra
\cite{Hoc77}, sometimes implicitly.  The late 1990s saw a resurgence
of explicit manifestations of Alexander duality, spurred in large part
by the Eagon--Reiner theorem \cite{ER}: possessing a linear free
resolution is dual to being Cohen--Macaulay.  In what follows, only
the most pertinent aspects of Alexander duality are mentioned; for
additional background---both mathematical and historical---see
\cite[Chapter~5]{cca}, including the Notes there.

The family of simplicial complexes in Definition~\ref{d:dual} was
identified by \Mustata\ \cite[Section~2]{Mus00} for the purpose of a
fundamental calculation, made independently by Terai \cite{Ter99}, of
the local cohomology of the polynomial ring~$S$ with support on a
(squarefree) monomial ideal.  The resulting Hilbert series formula
greatly resembles Hochster's formula for the local cohomology of a
Stanley--Reisner ring (Theorem~\ref{t:hilb}).  The relation between
the \Mustata--Terai formula and Hochster's formula is a functorial
effect of Alexander duality on local cohomology, known as local
duality with monomial support \cite[Section~6]{Mil00}.  This duality,
in turn, is a combinatorial special case of Greenlees--May duality
\cite{gm,MilGM}, which is generally an ajointness between the derived
functors of completion at an ideal~$I$ and of taking support on~$I$.

\begin{thm}\label{t:dual}
If $\aa \in \ZZ^n$ and $\1 = (1,\ldots,1)$, then the dual \v Cech
simplicial complex $\CD^{-\aa-\1}$ is Alexander dual to $\CD_\aa$
inside of the simplex on \mbox{$\{1,\ldots,n\} \minus \supp(\aa_-)$}.
\end{thm}
\begin{proof}
$\displaystyle\oF \in \CD^{-\aa-\1}
\begin{array}[t]{@{}l@{}}
       \iff \supp(\bb_-) = \oF\text{ for some }\bb\in -\aa-\1+\ol\Lambda_I
\\[1ex]\iff \bb_- = -\oF \text{ for some }\bb \in -\aa-\1+\ol\Lambda_I
\\[1ex]\iff -\oF \in \bb+\ZZ^F\text{ for some }\bb\in -\aa-\1+\ol\Lambda_I
\\[1ex]\iff -\oF \in -\aa-\1+\ol\Lambda_I+\ZZ^F
\\[1ex]\iff F-\1 \in -\aa-\1+\ol\Lambda_I+\ZZ^F
\\[1ex]\iff \aa + F \in \ol\Lambda_I + \ZZ^F
\\[1ex]\iff \aa \in \ol\Lambda_I + \ZZ^F
\\[1ex]\iff (IS_F)_\aa = \kk
\\[1ex]\iff (S_F/IS_F)_\aa = 0.
\\[1ex]\iff F - \supp(\aa_-) \not\in \CD_\aa.
\\[1ex]
\end{array}
$

\noindent
All of the equivalences are elementary; the only one demanding
verification of auxiliary data is the final line, for which it is
important to check that $F$ actually contains $\supp(\aa_-)$.  This is
not so bad, though: $\oF \in \CD^{-\aa-\1}$ implies that $\oF
\subseteq \supp((-\aa-\1)_-) = \supp((\aa+\1)_+)$, and this is the
complement of~$\supp(\aa_-)$.
\end{proof}

\begin{remark}\label{rk:Ext}
If \Mustata\ was trying to compute local cohomology with support on a
squarefree monomial ideal, how did he come upon a formula for
something equivalent to local cohomology---with maximal support---of
the quotient by an arbitrary monomial ideal?  Fairly easily: in
computing the limit that defines the local cohomology
$H^i_{I_\Delta}(S)$ supported on a squarefree monomial
ideal~$I_\Delta$, \Mustata\ was led to compute a simplicial formula
for $\Ext^i_S(S/I,S)$, where the monomial ideal~$I$ is a Frobenius
power of~$I_\Delta$.  But he also remarked that a similar formula
holds for all monomial ideals~$I$ \cite[Section~2]{Mus00}.  The
Alexander duality observed in this section is a direct result of
$\ZZ^n$-graded local duality (the ordinary kind, as opposed to with
monomial support):
\begin{align*}
  \Ext^i(S/I,S)_{-\aa-\1} &= \Ext^i(S/I,S(-\1))(\1)_{-\aa-\1}
\\                        &= \Ext^i(S/I,S(-\1))_{-\aa}
\\                        &= (H^{n-i}_\mm(S/I)_\aa)^\vee.
\end{align*}
\end{remark}

The formula for $\Ext^i(S/I,S)$ is perhaps the simplest formula in
this survey.

\begin{cor}\label{c:Ext}
If $\aa \in \ZZ^n$, then $\Ext^i(S/I,S)_\aa = \HH^{i-2}(\CD^\aa;\kk)$.
\end{cor}
\begin{proof}
The \emph{Alexander duality isomorphism} (see \cite[Theorem~5.6]{cca},
for example) states that, for two simplicial complexes $\Delta$ and
$\Delta^*$ that are Alexander dual inside of a simplex~$\sigma$, there
is a canonical isomorphism $\HH_{i-1}(\Delta;\kk) =
\HH^{|\sigma|-i-2}(\Delta^*;\kk)$.  Now apply the Alexander duality in
Theorem~\ref{t:dual} to the simplicial formula in
Theorem~\ref{t:cech}, using the local duality in Remark~\ref{rk:Ext}:
\begin{align*}
\Ext^{n-i}(S/I,S)_{-\aa-\1}
  &= (H^i_\mm(S/I)_\aa)^\vee
\\&= \HH_{i-|\supp(\aa_-)|-1}(\CD_\aa;\kk)
\\&= \HH^{n-i-2}(\CD^{-\aa-\1};\kk).
\end{align*}
Now, in the top and bottom lines, replace $n-i$ by~$i$ and $\aa$ by
$-\aa-\1$.  (Avoid trying to do this in the middle line, which has a
confusing $\supp(\aa_-)$ to deal with.)
\end{proof}

\begin{remark}
The formula in Corollary~\ref{c:Ext} seems to have been the first
simplicial local cohomology---or equivalently, Ext---formula for
arbitrary monomial ideals to appear in the literature
\cite[Section~2]{Mus00}, being a number of years earlier than the
next, which seems to have been \cite[Theorem~1]{takayama}.
\end{remark}

\section{Koszul simplicial complexes}\label{s:koszul}

As we have seen, verifying the Cohen-Macaulay condition by analyzing
the homology of a family of simplicial complexes, instead of a single
simplicial complex as in the squarefree case, can be accomplished
using a variety of simplicial complexes.  The simplest family,
however, is probably the following \cite[Definition~5.9]{cca}.

\begin{defn} \label{d:Kb}
Given a vector $\bb \in \NN^n$, set $\bb' = \bb - \supp(\bb)$.
For any monomial ideal~$I$, the \emph{(lower) Koszul simplicial
complex} of~$S/I$ in degree \mbox{$\bb \in \NN^n$}~is
\[
  \KD_\bb = \{\hbox{squarefree vectors } F \preceq \bb \mid
  \xx^{\bb'+F} \not\in I\}.
\]
Equivalently, if $\Lambda_I \subseteq \NN^n$ corresponds to the
standard monomials (Definition~\ref{d:std}),
\[
  \KD_\bb = \{F \preceq \supp(\bb) \mid \bb' + F \in \Lambda_I\}.
\]
\end{defn}

The motivation for the definition is analogous to the \v Cech
simplicial complexes: just as the algebraic chain complex of~$\CD_\aa$
is the graded piece of the \v Cech complex of~$S/I$ in degree~$\aa$,
the algebraic cochain complex of~$\KD_\bb$ is the graded piece of the
Koszul complex of~$S/I$ in degree~$\bb$.  In fact, the analogy between
these families of simplicial complexes is even closer.

\begin{lemma}\label{l:inclusion}
For $\bb \in \NN^n$ and $\aa = \bb-\1$, the Koszul simplicial complex
at~$\bb$ equals
\[
  \KD_\bb = \{F - \supp(\aa_-) \mid (S/I)_{\aa+F} = \kk\}.
\]
Consequently, $\CD_\aa$ is a subcomplex of\/~$\KD_\bb$.
\end{lemma}
\begin{proof}
The displayed formula for $\KD_\bb$ follows immediately from
Definition~\ref{d:Kb}, given that $(S/I)_{\aa+F}$ can only be nonzero
if $\aa+F$ has nonnegative coordinates, which occurs precisely when
$F$ contains~$\supp(\aa_-)$.  The containment $\CD_\aa \subseteq
\KD_\bb$ follows from Definition~\ref{d:cech}, given that
$(S_F/IS_F)_\aa$ is nonzero precisely when $(S/I)_{\aa+rF}$ is nonzero
for all $r \geq 1$.
\end{proof}

In certain special cases, as we shall see in Corollary~\ref{c:c-k},
below, the inclusion \mbox{$\CD_\aa \subseteq \KD_{\aa+\1}$} induces
an isomorphism on a particular cohomology group, and this cohomology
is precisely what is needed to detect the Cohen--Macaulay condition,
or its failure.  In general, a graded $S$-module~$M$ is
Cohen--Macaulay precisely when its minimal free resolution of~$M$ has
length equal to the codimension of~$M$ \cite[Theorem~13.37.2]{cca}.
Equivalently, the Betti numbers $\beta_{i,\aa}(M)$ must vanish in
homological degrees~$i$ greater than the codimension of~$M$, for
all~$\aa$.  Homologically, the Betti number in degree~$\aa$ and
homological degree~$i$ is the $\kk$-vector space dimension
\[
  \beta_{i,\aa}(M) = \dim_\kk \Tor^i_S(\kk,M)_\aa
\]
\cite[Lemma~1.32]{cca} (see also \cite[Definition~8.22]{cca}).
Therefore, the following generalization of Hochster's formula for
$\Tor^i_S(\kk,S/I_\Delta)_\aa$ to arbitrary monomial quotients~$S/I$
\cite[Theorem~5.11]{cca} yields the desired simplicial
characterization of the Cohen--Macaulay condition in
Corollary~\ref{c:koszul}.

\begin{thm}\label{t:koszul}
Given $\bb \in \NN^n$ with support $F = \supp(\bb)$, a simplicial
formula for the Betti number of~$S/I$ in degree~$\bb$ is given by the
$\kk$-vector space dimension of
\[
  \Tor_i(S/I,\kk)_\bb = \HH^{|F|-i-1} (\KD_\bb;\kk).
\]
\end{thm}
\begin{proof}
As mentioned after Definition~\ref{d:Kb}, this can be proved by direct
Koszul complex calculation.  (The proof of \cite[Theorem~5.11]{cca}
actually proceeds via Alexander duality; see Remark~\ref{rk:dual}.)
\end{proof}

\begin{cor}\label{c:koszul}
$S/I$ is Cohen--Macaulay if and only if $\HH^{|F|-i-1}(\KD_\bb;\kk) =
0$ for all $\bb \in \NN^n$ and all $i > n-d$, where $d = \dim(S/I)$
and $F = \supp(\bb)$.
\end{cor}

For the record, Hochster's Tor formula is the special case of
Theorem~\ref{t:koszul} in which $I = I_\Delta$ is squarefree.

\begin{cor}[Hochster's Tor formula]\label{t:tor}
The nonzero Betti numbers of~$S/I_\Delta$ lie only in squarefree
degrees~$F$, and we have
\[
  \Tor_i(S/I_\Delta,\kk)_F = \HH^{|F|-i-1}(\Delta|_F;\kk),
\]
where $\Delta_F$ is the subcomplex of all faces of~$\Delta$ whose
vertices lie in~$F$.
\end{cor}
\begin{proof}
This was originally proved in \cite{Hoc77}.  It follows immediately
from Theorem~\ref{t:koszul} and Definition~\ref{d:Kb}.  For additional
information and context, this result is also
\cite[Theorem~II.4.8]{Sta96}, \cite[Remark~5.5.5]{BH93}, and
\cite[Corollary~5.12]{cca}.
\end{proof}

\begin{remark}\label{rk:c-k}
For general vectors $\bb \in \NN^n$ and $\aa = \bb-\1$, the
relationship between the simplicial complexes $\CD_\aa$ and~$\KD_\bb$
in Lemma~\ref{l:inclusion} is mysterious, and perhaps difficult to
analyze (if such an analysis is possible, it would be interesting to
see).  This is to be expected, since the \v Cech simplicial complexes
compute local cohomology (Theorem~\ref{t:cech}), whereas the Koszul
simplicial complexes compute Tor into~$\kk$ (Theorem~\ref{t:koszul}),
which behaves quite differently from local cohomology, as a function
of the $\ZZ^n$-graded degree.
\end{remark}

Remark~\ref{rk:c-k} notwithstanding, the socle of the first nonzero
local cohomology precisely reflects the top Betti numbers.  This holds
quite generally for arbitrary finitely generated modules over regular
local or graded rings.  Here is the $\ZZ^n$-graded version.  Its proof
demonstrates that this result is a consequence of \emph{duality for
resolutions} \cite[Theorem~4.5]{Mil00} of $\ZZ^n$-graded modules.
That result was developed in the context of Alexander duality for
arbitrary (that is, not necessarily squarefree) monomial ideals and
modules, in response to the Eagon--Reiner theorem \cite{ER} and its
generalizations.  The point is that a minimal $\ZZ^n$-graded injective
resolution of a finitely generated module~$M$ contains data equivalent
to a minimal injective resolution of its Alexander dual, and therefore
also enough data to recover minimal free resolutions of both~$M$ and
its Alexander dual.  These considerations include the fact---key in
the proof of Theorem~\ref{t:c-k}---that a minimal free resolution
of~$M$ is equivalent to the maximal-support part of a minimal
injective resolution of~$M$.

\begin{thm}\label{t:c-k}
Let $M$ be finitely generated and $\ZZ^n$-graded of projective
dimension~$p$ ($\beta_{p,\bb}(M) \neq 0$ for some~$\bb \in \ZZ^n$, and
$\beta_{i,\bb}(M) = 0$ for $i > p$).  There is a natural~map
\[
  \Tor_i(M,\kk)_\bb \to H^{n-i}_\mm(M)_{\bb-\1}
\]
for all~$i$ and all~$\bb \in \ZZ^n$.  If $i = p$, and $\bb$ is any
degree such that the Betti number $\beta_{p,\bb}(M)$ is nonzero, then
the natural map is an isomorphism:
\[
  \Tor_p(M,\kk)_\bb \cong H^{n-i}_\mm(M)_{\bb-\1} \text{ if }
  \beta_{p,\bb} \neq 0.
\]
\end{thm}
\begin{proof}
Let $F_\spot$ be a minimal free resolution of~$M$.  Then
$\Tor_i(M,\kk)$ is the homology $H_i(F_\spot \otimes \kk)$.  In
constrast, by \cite[Theorem~4.5.5]{Mil00}, the $\ZZ^n$-graded
translate $H^{n-i}_\mm(M)(-\1)$ of the local cohomology up by~$\1$ is
the homology $H_i(F_\spot \otimes E_\kk)$, where $E_\kk$ is the
$\ZZ^n$-graded injective hull of~$\kk$, also known as the Matlis dual
of~$S$.  The existence of the natural map follows simply from the
canonical inclusion of~$\kk$ as the socle of~$E_\kk$.  The natural map
is an isomorphism when $i = p$ is maximal because, by minimality of
the resolution, the differentials of $F_\spot \otimes E_\kk$ are zero
on the socle.
\end{proof}

As a consequence, we find that the mysterious inclusion from
Lemma~\ref{l:inclusion} has at least one tractable feature.
Conveniently, this feature pertains precisely to the cohomology that
reflects the Cohen--Macaulay property for~$S/I$
(Corollary~\ref{c:koszul}).

\begin{cor}\label{c:c-k}
Let $S/I$ have projective dimension~$p$.  Assume that
\mbox{$\beta_{p,\bb}(S/I) \neq 0$}.  Writing $\aa = \bb-\1$ and $F =
\supp(\bb)$, the inclusion \mbox{$\CD_\aa \subseteq \KD_\bb$} induces
an isomorphism $\HH^{|F|-p-1} (\CD_\aa;\kk) =
\HH^{|F|-p-1}(\KD_\bb;\kk)$ on cohomology.\qed
\end{cor}

\begin{remark}
It is plausible that the \v Cech and lower Koszul simplicial complexes
in Corollary~\ref{c:c-k} could be equal, but neither a proof nor a
counterexample is known.  It is also possible that additional
hypotheses on $S/I$ might be required; starting with the case where
$S/I$ is Cohen--Macaulay is probably a good idea.
\end{remark}

\begin{remark}\label{rk:dual}
The formula in Theorem~\ref{t:koszul} comes in two flavors
\cite[Theorems~1.34 and~5.11]{cca} that are dual by definition, the
latter derived as an immediate consequence of the former, by Alexander
duality.  Hence it would be possible to present dual versions of the
results in this section, just as we did for the \v Cech simplicial
complexes~$\CD_\aa$ and their duals~$\CD^\aa$ in earlier sections.
However, the \v Cech simplicial complexes and their duals arose
separately, and their duality is neither by definition nor
self-evident, being a manifestation of local duality
(Remark~\ref{rk:Ext}).
\end{remark}

\raggedbottom

\end{document}